\documentclass{article}
\usepackage{mathrsfs,latexsym,amsmath}
\newcommand{\be}{\begin{equation}}
\newcommand{\ee}{\end{equation}}
\def\C{{\hbox{\bf C}}}

\def\Pr{{\hbox{\bf P}}}

\def\HH{{\cal H}}
\def\XX{{\mathscr X}_9}
\def\Q{{\hbox{\bf Q}}}
\def\Qbar{{\,\overline{\!\Q\!}\,}}

\def\Z{{\hbox{\bf Z}}}
\def\bmu{{\mbox{\boldmath$\mu$}}}
\def\vp{\varphi\0}
\def\ra{{\rightarrow}}
\def\0{^{\phantom0}}
\def\9{_{\phantom9}}
\def\m{_{\phantom-}}
\def\M{{\phantom{{M^|}^|}}}
\def\Gal{{\mathop{\hbox{\rm Gal}}\nolimits}}
\def\PSL{{\mathop{\hbox{\rm PSL}}\nolimits}}
\def\PGL{{\mathop{\hbox{\rm PGL}}\nolimits}}
\def\SL{{\mathop{\hbox{\rm SL}}\nolimits}}
\def\GL{{\mathop{\hbox{\rm GL}}\nolimits}}
\def\Aut{{\mathop{\hbox{\rm Aut}}\nolimits}}

\def\jinv{\hbox{$j$-invariant}}
\def\mat#1#2#3#4{{\genfrac(){0pt}{}{#1\;#2}{#3\;#4}}}
\begin{document}

\begin{Large}
\centerline{Elliptic curves with $3$-adic Galois representation}
\centerline{surjective mod~$3$ but not mod~$9$}
\end{Large}
\vspace*{2ex}

\centerline{Noam D. Elkies\footnote{
  Supported in part by NSF grants DMS-0200687 and DMS-0501029.
  }
}

\begin{quote}
\begin{small}
{\bf Abstract.}  Let $E$\/ be an elliptic curve over~$\Q$,
and $\rho_l: \Gal(\Qbar/\Q) \ra \GL_2(\Z_l)$
its \hbox{$l$-adic} Galois representation.
Serre observed that for $l\geq 5$ there is no proper closed subgroup
of $\SL_2(\Z_l)$ that maps surjectively onto $\SL_2(\Z/l\Z)$,
and concluded that if $\rho_l$ is surjective mod~$l$\/
then it is surjective onto $\GL_2(\Z_l)$.
We show that this no longer holds for $l=3$
by describing a modular curve~$\XX$ of genus~$0$
parametrizing elliptic curves for which $\rho_3$
is not surjective mod~$9$ but generically surjective mod~$3$.
The curve $\XX$ is defined over~$\Q$, and the modular cover
$\XX \ra X(1)$ has degree~$27$, so $\XX$ is rational because $27$ is odd.
We exhibit an explicit rational function~$f\in\Q(x)$ of degree~$27$
that realizes this cover.  We show that for every $x\in\Pr^1(\Q)$,
other than the two rational solutions of $f(x)=0$, the elliptic curves
with \jinv\ $f(x)$ have $\rho_3$ surjective mod~$3$ but not mod~$9$.
We determine all nonzero integral values of~$f(x)$, and exhibit
several elliptic curves satisfying our condition on~$\rho_3$,
of which the simplest are the curves
$Y^2=X^3-27X-42$ and $Y^2+Y=X^3-135X-604$
of conductors $1944 = 2^3 3^5$ and $6075 = 3^5 5^2$ respectively.
\end{small}
\end{quote}

{\bf 0.~Introduction.}
Let $E$\/ be an elliptic curve over~$\Q$,
and $\rho_l: \Gal(\Qbar/\Q) \ra \GL_2(\Z_l)$
its \hbox{$l$-adic} Galois representation.
Serre observed in 1968 \cite[IV, 3.4, Lemma~3]{Serre}
that for $l\geq 5$ there is no proper closed subgroup
of $\SL_2(\Z_l)$ that maps surjectively onto $\SL_2(\Z/l\Z)$,
and concluded that if $\rho_l$ is surjective mod~$l$\/
then it is surjective onto $\GL_2(\Z_l)$.
He noted \cite[IV, 3.4, Exercise~3]{Serre}
that for $l=3$ there exists a subgroup $G \subset \SL_2(\Z/9\Z)$
such that the restriction to~$G$\/ of the \hbox{reduction-mod-3} map
$\SL_2(\Z/9\Z) \ra \SL_2(\Z/3\Z)$ is an isomorphism.
The preimage of~$G$ in $\SL_2(\Z_3)$ is then a proper closed subgroup
that maps surjectively to $\SL_2(\Z/3\Z)$.
This suggests that there could be curves~$E$\/
for which $G$\/ is the image of $\rho_3 \bmod 9$,
making $\rho_3$ surjective mod~$3$ but not mod~$9$.
Serre does not raise this question explicitly,
and it does not seem to have been addressed elsewhere in the literature;
I thank Grigor Grigorov for drawing my attention to it.
In this paper we answer the question by showing that there exist
infinitely many $j \in \Q$ for which an elliptic curve
of \jinv~$j$ must have $\rho_3$ surjective mod~$3$ but not mod~$9$.
The simplest examples are $j=4374$, $j=419904$, and $j=-44789760$.
In general $j$ is the value of a rational function $f(x)$
of degree~$27$ at all but finitely many $x\in\Pr^1(\Q)$.

Such curves $E$\/ are parametrized by a modular curve
$\XX = X(9) / G$.  The natural cover $\XX \ra X(1)$ has degree~$27$,
and our rational function~$f$\/ arises as the pullback to~$\XX$
of the \hbox{degree-$1$} function~$j$ on~$X(1)$.
It is easy to check from the Riemann-Hurwitz formula
that $\XX$ has genus~zero.  The challenge is to prove
that $\XX$ is defined over~$\Q$ and to compute~$f(x)$
for some choice of rational coordinate~$x$ on~$\XX$.
(Once $\XX$ is known to be defined over~$\Q$,
it is automatically isomorphic with~$\Pr^1$ over~$\Q$,
becoause it supports the rational function~$j$ of odd degree.)
We prove the rationality in section~1,
and compute~$x$ using products of Siegel functions in section~2.
Such products are modular units on~$X(N)$ with known
\hbox{$q$-expansions}; since the cusps of~$\XX$ are not rational
we must also find a fractional linear transformation over~$\Qbar$
that takes the modular unit to a function~$x$ defined over~$\Q$.
In section~3 we discuss elliptic curves~$E$\/
with \jinv s obtained by specializing~$f$,
and use explicit computation of curves $X(9)/H$\/ with $H \subset G$\/
to show that all such curves with $x\in\Pr^1(\Q)$
satisfy our condition on~$\rho_3$ except for those with $j=0$.
Finally in section~4 we determine the finite set $f(\Pr^1(\Q)) \cap \Z$,
and exhibit some specific elliptic curves~$E$\/
whose \jinv s are these integral values.

\vspace*{2ex}

{\bf 1.~The group~$G$\/ and the curve~$\XX$.}
The group $\SL_2(\Z/3\Z)$ is generated by
$$
S = \left( \!
  {\begin{array}{rr} 0 & \!\! 1 \cr -1 & \!\! 0 \end{array}}
  \right),
\qquad
T = \left(
  {\begin{array}{rr} 1 & \! 1 \cr 0 & \! 1 \end{array}}
  \right),
$$
the images mod~$3$ of the standard generators
of $\SL_2(\Z)$.  These generators of~$\SL_2(\Z/3\Z)$ satisfy
$$
S^2 = (ST)^3 = -I, \quad T^3 = I.
$$
To lift $\SL_2(\Z/3\Z)$ to a subgroup~$G$\/ of~$\SL_2(\Z/9\Z)$,
it is enough to lift $S,T$\/ to matrices mod~$9$
satisfying the same relations.\footnote{
  Warning: This approach to lifting $\SL_2(\Z/n\Z)$
  works only for $n\leq 5$ (including $n=4$),
  where $\langle s,t \mid s^2 = (st)^3 = t^n = 1 \rangle$
  is a presentation of $\PSL_2(\Z/n\Z)$.
  For $n>5$ more relations must be checked.
  }
A direct search finds $27$ such lifts $\widetilde S, \widetilde T$,
all equivalent under conjugation in $\SL_2(\Z/9\Z)$.
These yield $27$ choices of lift of $\PSL_2(\Z/3\Z)$
to a subgroup $G / \{\pm 1\}$\/ of $\PSL_2(\Z/9\Z) = \Aut(X(9))$,
all conjugate in $\Aut(X(9))$.
Hence the quotients of $X(9)$ by these subgroups
are all equivalent under $\Aut(X(9))$.  We choose
$$
\widetilde S = \left(
  {\begin{array}{rr} 0 & \!\! 2 \cr 4 & \!\! 0 \end{array}}
  \right),
\qquad
\widetilde T = \left( \!
  {\begin{array}{rr} 4 & \! 1 \cr -3 & \! 4 \end{array}}
  \right),
$$
and set
$G = \langle \widetilde S, \widetilde T \kern.08ex\rangle
\subset \SL_2(\Z/9\Z)$
and $\XX = X(9) \bigl/ (G / \{\pm1\}) \bigr.$.

We can then calculate the genus of~$\XX$ by applying
the Riemann-Hurwitz formula either to the quotient map $X(9)\ra\XX$,
using the fact that $X(9)$ has genus~$10$,
or to the covering map $\XX \ra X(1)$.
In the quotient map, each of the three involutions in~$G$\/
has six fixed points, and each of the four \hbox{$3$-element} subgroups
has three; thus the ramification divisor has degree
$3 \cdot 6 + 4 \cdot 6 = 42 = 2|G| + 2(10-1)$,
whence $\XX$ has genus~$0$.  The covering map has degree~$27$
and is unramified except above the cusp $j=\infty$
and the elliptic points $j=0$, $j=1728$.
We find that these points have preimages with multiplicities
$9^3$, $3^8 1^3$, and $12^2 1^3$ respectively
(using $m^c$ as a standard shorthand for $c$ preimages
of multiplicity~$m$).  Hence the ramification divisor
has degree $3 \cdot 8 + 8 \cdot 2 + 12 = 52 = 2(27-1)$,
so again we conclude that $\XX$ has genus~$0$.

In particular, each of the cusps of $X(9)$
has trivial stabilizer in~$G$.  To explain this,
note that the stabilizer in $\PSL_2(\Z/9\Z)$ of each cusp of~$X(9)$
is conjugate to the group of matrices $\{\pm\mat1a01\}$;
if $G$\/ had a nontrivial intersection with this group
then $G$\/ would contain $\{\pm\mat1301\}$,
contradicting the requirement that the reduction map
$G \ra \PSL_2(\Z/3\Z)$ be bijective.
It follows that $\XX$ has $27/9=3$ cusps.

To obtain a model of~$\XX$ defined over~$\Q$, we must extend~$G$\/
to a group $G' \subset \GL_2(\Z/9\Z)$ such that $G' \rhd G$\/ and
the determinant map $G'/G \ra (\Z/9\Z)^*$ is an isomorphism.
For the preimages of the squares in~$(\Z/9\Z)^*$
we use the invertible multiples of the identity
(recall that $-I$\/ is already in~$G$\/).
It remains to lift $\pm\mat{-1}0{\m0}1$ from $\GL_2(\Z/3\Z)$
to $\GL_2(\Z/9\Z)$, and we calculate that the unique choice
that normalizes~$G$\/ is $\pm\mat{-1}0{\m0}1$.
Together with the invertible scalar matrices
this yields a unique extension~$G'$ satisfying our conditions,
and thus a unique structure for $\XX$ as a modular curve over~$\Q$.

In general, even if a curve of genus~$0$ is known to be defined
over~$\Q$ it need not be isomorphic with $\Pr^1(\Q)$.
But in our case the curve supports \hbox{$\Q$-rational}
divisors of odd degree (such as the preimage of the cusp,
or indeed the preimage of any rational point on $X(1)$),
so $\XX$ must be isomorphic with $\Pr^1(\Q)$.
In the next section we choose an isomorphism, and compute $j$
as a rational function of a rational coordinate on~$\XX$.

While $\XX$ is defined over~$\Q$, its three cusps are not rational.
More precisely, they are conjugate over the cyclic cubic extension
$K := \Q(\zeta+\zeta^{-1})$, where $\zeta = e^{2 \pi i/9}$,
a primitive 9th root of unity.  To see this, we check that
$G'$ contains an element $\pm\mat{a}{b}01$ for $a=\pm1$
(when $b=0$) but no other $a\in(\Z/9\Z)^*$.
(Since $G'$ contains $\mat{-1}0{\m0}1$,
it is enough to check this for $a=4$ and $a=7$,
and then we can multiply by the scalar matrix $a^{-1} \in G'$
to reduce to the corresponding statement for~$G$.)
Thus the subgroup of $\Gal(\Q(\zeta)/\Q)$ that fixes the cusp $i\infty$
consists of the images of $\pm1$ under the standard identification
of $(\Z/9\Z)^*$ with $\Gal(\Q(\zeta)/\Q)$.  Hence that cusp
is defined over~$K$\/ but not over~$\Q$, and its conjugates
under $\Gal(K/\Q)$ must be the other two cusps of~$\XX$.

The fact that the cusps of~$\XX$\/ are not rational will make
our computation of rational functions on~$\XX$\/ somewhat trickier:
we can still expand these functions in powers of~$q$,
but the coefficients will in general be contained only in~$K$,
not in~$\Q$, even for a function in $\Q(\XX)$.

\vspace*{2ex}

{\bf 2.~Computing $x$ and $f(x)$.}
For $\tau$ in the upper half-plane~$\HH$,
let $q = e^{2 \pi i \tau}$ as usual,
and let $q_9 = e^{2 \pi i \tau/9}$,
a local parameter for $X(9)$ and~$\XX$ at the cusp $\tau = i\infty$,
with $q = q_9^9$.
The curve $X(3)$ has genus zero with rational coordinate
$$
H_3 = \left(\frac{\eta(\tau/3)}{\eta(3\tau)}\right)^{\!3} + 3
= q_9^{-3} (1 + 5 q - 7 q^2 + 3 q^3 + 15 q^4 - 32 q^5 + + - \cdots);
$$
the curve $X(9)$ is a $(\Z/3\Z)^3$ cover of~$X(3)$
whose function field is obtained from $\C(H_3)$ by adjoining
cube roots of $H-3$, $\zeta^3 H - 3$, and $\zeta^{-3} H - 3$.
The group $\PSL_2(\Z/3\Z)$ acts on $X(3)$
by the fractional linear transformations of~$H$\/
that permute $\{ \infty, 3, 3\zeta^3, 3\zeta^{-3} \}$, which are
the transformations that preserve $j = H^3 (H^3+216) / (H^3-27)^3$.
(Recall that we set $\zeta = e^{2 \pi i/9}$,
so $1$ and $\zeta^{\pm 3}$ are the cube roots of unity.)
To compute $\XX$, we may lift these transformations to elements of~$G$\/
in the group $\PSL_2(\Z/9\Z)$ of automorphisms of $X(9)$
and look for a rational function~$x$ of degree~$1$
on the quotient $\XX = X(9) / G$.
We could try to find~$x$
using our above description of the function field of $X(9)$
as a $(\Z/3\Z)^3$ extension of $\C(H_3)$,
but this seemed like an unpleasant project. 
Instead we use modular units on~$X(9)$, that is,
modular functions whose divisors are supported on the cusps.

For integers $a,b$ we have a function $s_{a,b}$ on~$\HH$\/ defined by
$$
s(a,b) = q_9^\alpha (1 - \zeta^b q_9^a)
\prod_{n=1}^\infty (1 - \zeta^b q_9^{9n+a}) (1 - \zeta^{-b} q_9^{9n-a}),
$$
where $\alpha = (a^2/9 + a - 9/6) / 2 = 18 B_2(a/9)$.
This is a modular unit, namely
the ``Siegel function'' of~\cite[p.29]{KL} with parameters $(a/9, b/9)$,
multiplied by a power of $-\zeta$ that is irrelevant for our purpose.
Taking $N=9$ in~\cite[p.68, Thm.~4.1]{KL},
we see that a product $F = \prod_r s(a_r,b_r)^{m_r}$ of such functions
is a modular function on~$X(9)$ if and only if
the $(a_r,b_r)$ and $m_r$ satisfy the ``quadratic relations''
$$
\sum_r m_r\0 a_r^2 \equiv \sum_r m_r\0 b_r^2 \equiv \sum_r m_r a_r b_r
\equiv 0 \bmod 9.
$$
Note that this condition depends only on the $(a_r,b_r) \bmod 9$,
and is invariant under changing $(a_r,b_r)$
to $(-a_r,-b_r)$ for some~$r$;
Indeed $s(-a,-b)$ is a scalar multiple of $s(a,b)$,
as is $s(a',b')$ if $(a',b') \equiv (a,b) \bmod 9$.

The pairs $\{(a,b), (-a,-b)\} \in (\Z/9\Z)^2$
with at least one of $a,b$\/ not divisible by~$3$
correspond bijectively with the cusps of $X(9)$:
the cusps are orbits of $\Gamma(9)$ acting on $\Pr^1(\Q)$,
and if $\gcd(a,b)=\gcd(a',b')=1$ for some integers $a,b,a',b'$
then $a/b$ and $a'/b'$ are in the same orbit
if and only if $(a',b') \equiv (a,b)$ or $(-a,-b) \bmod 9$.
This labeling is consistent with the action of $\PSL_2(\Z/9\Z)$
on cusps, and dual to its action on the modular units $s(a,b)$
(up to $\bmu_{18}$ factors).  We choose one of the three orbits
of the action of $G/\{\pm 1\}$ on the cusps of $X(9)$,
and take $(a_r,b_r)$ corresponding to the twelve cusps in the orbit.
We check that these satisfy the quadratic relations with $m_r=1$,
so the product~$F$\/ of the twelve functions $S(a_r,b_r)$ 
is a modular function on~$X(9)$ whose divisor is invariant under~$G$.

There is then a homomorphism $\chi: G \ra \C^*$ such that
$F(g(\tau)) = \chi(g) F(\tau)$ for all $g\in G$\/ and $\tau \in \HH$.
We claim that $\chi$ is trivial.
Indeed, if it were nontrivial we could take for~$g$
a \hbox{$3$-cycle} in~$G$,
and for~$\tau$ a preimage of a fixed point on $X(9)$ of~$g$,
and conclude that $\tau$ is a zero or pole of~$F$,
which is impossible because $\tau$ is not a cusp.
Therefore $\chi$ is trivial as claimed, whence $F$\/ descends
to a function on $\XX$ whose only poles or zeros are at the cusps.

We find that $F$\/ has a simple pole at one cusp,
a simple zero at another, and neither zero nor pole at the third.
In particular, $F$\/ is a function of degree~$1$ on~$\XX$.
Choosing the orbit of $(a,b)=(1,0)$, corresponding to the cusp $i\infty$,
we find that $F$\/ has a pole at that cusp,
and calculate the \hbox{$q$-expansion}
$$
F = q_9^{-1} - 1 + c_4\0 q_9\0 + c_1\0 q_9^2
+ (c_2\0 + 2) q_9^3 + c_4\0 q_9^4 \cdots,
$$
where $c_m := \zeta^m + \zeta^{-m}$, a unit in~$K$.
The \hbox{$q$-expansions} of the products corresponding
to the other two orbits then let us recognize those products
as fractional linear transformations of~$F$,
namely $1/(F-c_2+1)$ and \hbox{$(1-c_2)/F$}.

It follows that $F$\/ takes the values $0$ and $1-c_2$
on the other two cusps of~$\XX$.  Therefore $j$,
considered as a rational function of~$F$,
must have its poles at $F=0,1-c_2,\infty$,
each with multiplicity~$9$.  Indeed we compute that
$F^9 (F - c_2 + 1)^9  j$ is a polynomial of degree~$27$ in~$F$\/
to the accuracy allowed by our \hbox{$q$-expansions}
(which extend far enough beyond the constant term of that polynomial
to provide a sanity check on our computations).

It remains to find a fractional linear transformation
with coefficients in~$K$\/ that, when applied to~$F$,
yields a coordinate~$x$ on~$\XX$ such that $j \in \Q(x)$.
Thus $x$ must map the three cusps to a $\Gal(K/\Q)$ orbit in~$K$.
We may choose any ordered orbit, and then $x$ is determined uniquely,
because $\PGL_2$ acts simply \hbox{$3$-transitively} on~$\Pr^1$.
(The order must be consistent with the Galois action on the cusps.)
We then apply a $\PGL_2(\Q)$ transformation so that the map
$\XX \ra X(1)$ is represented by a rational function
with small coefficients.  This leads us to
\begin{eqnarray*}
x & \!\! = \!\! & \frac{-c_1 F + 1 - c_2}{F - c_1 + 3(1-c_2)}
\\
& \!\! = \!\! &
-c_1 + (2c_2+c_4) q_9\0 + 3(1-c_1) q_9^2 + (6-7c_1+c_2) q_9^3
+ (15-16c_1+7c_2) q_9^4 \cdots,
\end{eqnarray*}
when $j=f(x)$ with
\begin{eqnarray*}
f(x) & \!\! = \!\! & - \frac
  {3^7 (x^2-1)^3 (x^6+3x^5+6x^4+x^3-3x^2+12x+16)^3 (2x^3+3x^2-3x-5)}
  {(x^3-3x-1)^9}
\\
  & \!\! = \!\! & 1728 - \frac
  {3^3 A^2(x) B^2(x) (2x^3-3x^2+4)} {(x^3-3x-1)^9},
\end{eqnarray*}
where $A(x)$, $B(x)$ are the sextic polynomials
\begin{eqnarray*}
A(x) & \!\! = \!\! & x^6+6x^5+4x^3+12x^2-18x-23,
\\
B(x) & \!\! = \!\! & 7x^6+24x^5+18x^4-26x^3-33x^2+18x+28.
\end{eqnarray*}

\vspace*{2ex}

{\bf 3.~The elliptic curves parametrized by~$\XX$.}
Now let $E/\Q$ be an elliptic curve with \jinv\ $f(x)$
for some $x\in\Pr^1(\Q)$.
Assume that its \hbox{$3$-adic} Galois representation $\rho_3$
is surjective mod~$3$.  Then its image mod~$9$ is (a conjugate of)
the proper subgroup $G'$ of~$\GL_2(\Z/9\Z)$, because this image
is contained in~$G$\/ and its determinant maps surjectively to $\Z_3^*$.
In particular, $\rho_3$ is surjective mod~$3$ but not mod~$9$,
as desired. 

We claim that the \hbox{mod-$3$} condition on~$\rho_3$
is satisfied by except at $x=\pm1$, the points at which $f(x)=0$,
the \jinv\ of an elliptic curve with complex multiplication.

For any elliptic curve~$E/\Q$,
the representation $\rho_3$ is surjective mod~$3$ if and only if
the intersection of its image in $\PGL_2(\Z/3\Z)$ with $\PSL_2(\Z/3\Z)$
is contained in the \hbox{$4$-element} normal subgroup
or in a cyclic subgroup of order~$3$.
In the former case, the \jinv~$j_E$ of~$E$\/ is a cube,
and this necessary condition is also sufficient unless $j_E=0$.
If $f(x) \neq 0$ then $f(x)$ is a cube if and only if
$3(2x^3+3x^2-3x-5) = z^3$ for some $z\in\Q$.
This curve of genus~$1$ has no rational points
due to a \hbox{$3$-adic} obstruction: we have
$$
2x^3+3x^2-3x-5 = 3 (x+1)^3 - (x+2)^3,
$$
so the \hbox{$3$-adic} valuation of $3(2x^3+3x^2-3x-5)$
is never divisible by~$3$.  Thus $f(x)$ is never a nonzero cube.
The latter case holds if and only if
$E$\/ admits a rational \hbox{$3$-isogeny}.
Such $x$ are parametrized by a curve
$\XX' = \XX / \langle \widetilde T \rangle$,
the quotient of $X(9)$ by a \hbox{$3$-cycle} in~$G$.
This curve has genus~$3$,
so there are only finitely many such~$x$ by Mordell-Faltings.
For a general curve of genus~$3$ it is not known how to provably list
all the rational points.  But here we are lucky:
we can give a rational map of degree~$3$ from~$\XX'$
to the elliptic curve $Y^2+Y=X^3$, which is known to have rank zero.
Pulling back each of this curve's three rational points to~$\XX'$
then yields a set of \hbox{$\Qbar$-rational} points
that must contain all the \hbox{$\Q$-rational} ones.
To complete the proof of our claim we shall observe that there are
only two \hbox{$\Q$-rational} points, one for each of $x=\pm1$.

A simple model for $\XX'$ is
\be
3 z_2^3 = z_1\0 (z_1^3 + 3 z_1^2 - 6 z_1\0 + 1),
\label{eq:XX'}
\ee
and the map to $Y^2+Y=X^3$ can be given by
\be
(X:Y:1) = (z_1 (z_1+1) z_2 : 3 z_1^2 : z_2^3).
\label{eq:XX'map}
\ee
The rational points on~$\XX'$ are $(0,0)$ and the point at infinity.
The equation (\ref{eq:XX'}) for~$\XX'$ shows that this curve
has an automorphism~$\sigma$ of order~$3$
that fixes~$z_1$ and multiplies $z_2$ by a cube root of unity.
This automorphism arises from the element
$3\widetilde T - 2 = \mat1301$ of $\PSL_2(\Z/9\Z)$,
which commutes with~$T$\/ and thus descends from $\Aut(X(9))$
to an automorphism of~$\XX'$.
The quotient of~$\XX'$ by~$\langle\sigma\rangle$
is the \hbox{$z_1$-line}, which covers the \hbox{$j$-line $X(1)$}
with degree~$36$ and the curve $X_0(3)$ with degree~$9$.
The latter map can be realized by the rational function
$
27 \bigl(
  (z_1^3 + 3 z_1^2 - 6 z_1\0 + 1) / (z_1^3 - 6 z_1^2 + 3 z_1\0 + 1)
\bigr)^{\!3}.
$

We next outline the computation of the model~(\ref{eq:XX'})
and the map~(\ref{eq:XX'map}).  We begin with the modular units
\be
s(1,0) s(4,6) s(4,3),
\quad
s(4,0) s(7,6) s(7,3),
\quad
s(7,0) s(1,6) s(1,3).
\label{eq:sss}
\ee
The first of these corresponds to the orbit of the cusp $i\infty$
under the \hbox{$3$-cycle}~$\widetilde T$, and the others are obtained
by multiplying $\{(1,0), (4,6), (4,3)\}$ by~$4$ and $7 \bmod 9$.
We calculate that these products do not satisfy the quadratic relations,
but the quotient of any two of them does, giving a rational function 
of degree~$3$ on~$\XX'$.  Using the \hbox{$q$-expansions} we find that
the three functions in~(\ref{eq:sss}) are linearly dependent,
and thus that their pairwise quotients are all related by fractional
linear transformations over~$K$.  As we did for~$x$,
we find two linear combinations of the functions in~(\ref{eq:sss})
whose quotient is a \hbox{degree-$3$} function on~$\XX'$
defined over~$\Q$.  We call this function
$$
y = - c_2 + (c_4-c_2+3) q_9\0 + (3c_4-6c_7+9) q_9^2
+ (10c_4-22c_7+27) q_9^3 \cdots,
$$
and use the $q$-expansions to find the coefficients of
a polynomial identity $P(x,y)=0$ of bidegree $(3,4)$ in~$(x,y)$.

This is a model for~$\XX'$ in $\Pr^1 \times \Pr^1$, but it is not smooth.
We find that it has three double points, at $x = y = -c_1, -c_2, -c_4$.
Thus the holomorphic differentials on~$\XX'$ are the forms
$Q(x,y) \, dx/P_y$ for $Q$\/ in the \hbox{$3$-dimensional} space
of polynomials of bidegree~$(1,2)$
that vanish at the three singularities.
We interpret these as modular cuspforms of weight~$2$ on~$\XX'$
by writing $q \, dx/dq$ in place of~$dx$,
and find a basis for the space of such cuspforms:
\begin{eqnarray*}
\vp_1 & \!\! = \!\! &
c_4 q - 3 q^3 - 2 c_2 q^4 - c_1 q^7 + 6 q^{12} + 5 c_2 q^{13}
+ 4 c_1 q^{16} - 7 c_4 q^{19} + 3 q^{21} \cdots,
\\
\vp_2 & \!\! = \!\! &
q + (c_4 - c_1) q^2 + q^4 + (2c_2 + c_4) q^5 + 2q^7 + (c_2+c_4) q^8
- 3 q^{10} \cdots,
\\
\vp_3 & \!\! = \!\! &
(c_4 - c_2) q - 3 q^2 + (c_2 - c_1) q^4 + 3 q^5 + (2c_1-2c_4) q^7
+ 3 q^8 + (3c_2-3c_4) q^{10} \cdots.
\end{eqnarray*}
The affine model (\ref{eq:XX'}) is then obtained by taking
$(z_1:z_2:1) = (\vp_2 + \vp_1 : \vp_3 : 2\vp_2 - \vp_1)$.
The map~(\ref{eq:XX'map}) was obtained by integrating
the CM~form~$\vp_1$; it can also be seen in~(\ref{eq:XX'})
by writing the cubic factor as $(z_1+1)^3 - 9 z_1$.

\vspace*{2ex}

{\bf 4.~Numerical examples.}
Besides $f(\pm 1) = 0$, there are seven other integers
obtained by evaluating $f(x)$ at points of~$\Pr^1(\Q)$
of small height.  For each of those, we list $x$, $f(x)$,
and one of the elliptic curves $E$\/ of \jinv\ $f(x)$
and minimal conductor~$N$\/:

\vspace*{2ex}

\centerline{
$
\begin{array}{c|c|c|c}
  x   &       j=f(x)     &           E              &    N
 \\ \hline
 1/0  &    \M  4374  \M  &       [0,0,0,-27,-42]    & 2^3 3^5
 \\
 -2\m &       419904     &      [0,0,0,-162,792]    & 2^8 3^5
 \\
  0   &    -44789760\m   &      [0,0,1,-135,-604]   & 3^5 5^2
\\
-1/2\m&   15786448344    &   [0,0,0,-5427,153882]   & 2^5 3^5
\\
  2   &  24992518538304  & [0,0,0,-201042,34695912] & 2^8 3^5 17^2
\\
-3/2\m& -92515041526500\m& [0,0,0,-1126035,459913278] & 2^3 3^5 19^2
\\
-1/3\m& -70043919611288518656\m &  [0,0,1,-1127379978,-14569799990728] &
                                                    3^5 97^2 101^2
\end{array}
$
}

\vspace*{1ex}

The smallest conductor here is $1944 = 2^3 3^5$, still too large
to appear in Cremona's published tables~\cite{Cremona}
of curves of conductor $\leq 1000$.  But Cremona
has pursued his computations up to $10^5$ and beyond
(see \cite{Cremona1} for the status as of mid-2006),
enough to find our first curve
as well as those of conductors $6075 = 3^5 5^2$, $7776 = 2^5 3^5$,
and $62208 = 2^8 3^5$.
Curves with $j=4374$, $j=419904$, and $j=15786448344$
already appeared in the tables of elliptic curves
with good reduction away from $2$~and~$3$,
compiled in~1966 by F.B.~Coghlan
and published as ``Table 4'' in \cite[p.123]{Antwerp}:
see rows 52, 84, and~86 of ``Table~4a'' \cite[p.125]{Antwerp},
and \cite[p.75]{Antwerp} for the attribution to Coghlan.

A search up to height $256$ found no more integral values of $f(x)$.
Since $f$\/ has three distinct poles, there can be only finitely many~$x$
for which $f(x) \in \Z$.  We claim that in fact we have found them all.
Suppose $x=m/n$ in lowest terms.  The resultant of the numerator
and denominator of $f(x)$ is $3^{486}$, so when we write $f(m/n)$
as the quotient of homogeneous polynomials in $m,n$ the denominator
$(m^3 - 3mn^2 - n^3)^9$ must be $\pm 3^A$ for some~$A$.
Thus $m^3 - 3mn^2 - n^3 = \pm 1$ or~$\pm3$, because
the only $(m,n)\in\Z^2$ for which $9 \mid m^3 - 3mn^2 - n^3$
are those for which $3|m$ and $3|n$.
The only cases of $m^3 - 3mn^2 - n^3 = \pm 3$ are $x=1,-2,-1/2$,
because these are the only rational points on the elliptic curve
$x^3-3x-1=3z^3$ (isogenous with the cubic Fermat curve).
For $m^3 - 3mn^2 - n^3 = \pm 1$ one must work harder
because the elliptic curve $x^3-3x-1=z^3$ has rank~$1$.
Fortunately this work was already done by Ljunggren,
who proved in 1942 that the Thue equation $m^3 - 3mn^2 - n^3 = 1$
has only the six solutions $(1,0)$, $(0,-1)$, $(-1,1)$,
$(2,1)$, $(1,-3)$, $(-3,2)$, corresponding to values of~$x$
already listed above.  (See \cite[\S2]{Lj},
cited by Nagell~\cite{Nagell} who also notes the connection
between $X^3+Y^3=Z^3$ and $m^3 - 3mn^2 - n^3 = 3$.)
The solutions of $m^3 - 3mn^2 - n^3 = -1$ are obtained from these six by
changing each $(m,n)$ to $(-m,-n)$, which yields the same values of $x=m/n$.
Thus our list contains all the nonzero integral \jinv s
of elliptic curves parametrized by~$\XX$.

\end{document}